\documentclass{article}

\usepackage{amsfonts,amsmath}

\def\Z{\mathbb {Z}}
\def\qed{{\hfill\vrule height 5pt width 5pt depth 0pt}}
\def\phi{\varphi}
\newtheorem{lemma}{Lemma}
\newtheorem{theorem}[lemma]{Theorem}

\title{Group Actions in Number Theory}

\author{Benjamin V. Holt 
\\Department of Mathematics
\\Humboldt State University
\\Arcata, CA 95521 USA
\\E-mail: bvh6@humboldt.edu
\and         
Tyler J. Evans 
\\Department of Mathematics
\\Humboldt State University
\\Arcata, CA 95521 USA
\\E-mail: te8@humboldt.edu
}

\begin{document}

\maketitle

\section{Introduction}

Students having had a semester course in abstract algebra are exposed to the elegant way in which finite group theory leads to proofs of familiar facts in elementary number theory.  In this note we
offer two examples of such group theoretical proofs using the action of a group on a set.  The first is Fermat's little theorem and the second concerns a well known identity involving the famous Euler phi function. The tools that we use to establish both results are sometimes seen in a second semester algebra course in which group actions are studied.  Specifically, we will use the class equation of a group action and Burnside's theorem.

\section{Fermat's Little Theorem}

A well known consequence of the class equation of a group action asserts
that if $G$ is a $p$-group (that is, $G$ is a finite group of order $p^n$
for some integer $n\ge 1$ and a prime integer $p$), and $G$ acts on a
finite set $S$, then the number of elements in $S$ is congruent to the
number of fixed points of the action modulo $p$.  Recall that an element
$s\in S$ is a fixed point of the action if $gs=s$ for all $g\in G$.  The
set of fixed points is usually denoted by $S^G$, and so using this
notation, the aforementioned theorem asserts that 
\begin{equation}
\label{fact} |S|\equiv |S^G|\pmod p. 
\end{equation}
This seemingly obscure result appears to have limitless utility in group
theory!  A line of argumentation due to R.~J.~Nunke (c.f~\cite{Fraleigh,Hungerford}) uses
(\ref{fact}) repeatedly to establish the three famous Sylow theorems in
elementary group theory.  In
this section, we use (\ref{fact}) to obtain a new proof of the following
well known result in number theory.

\begin{theorem}[Fermat's little theorem]\label{fermat}
If $a\ge 1$ is any integer and $p$ is a prime, then $a^p\equiv a \pmod p$.
\end{theorem}

{\bf Proof.} Let $a\ge 1$ be an integer and let $A=\{1,\cdots,a\}$.  Let
$S=A^p$ and let a cyclic group $G$ of order $p$ act on
$S$ by cyclic permutation of the entries of an element in $S$.  Easily, this is
a well defined action and $(a_1,\cdots,a_p)\in S$ is a fixed point if
and only if $a_1=\cdots = a_p$.  Therefore $|S^G|=a$.  Since
$|S|=a^p$, an application of (\ref{fact}) completes the proof. \qed

Of course, Fermat's little theorem holds for all integers $a$, but the construction in the our proof is not valid for $a\le 0$.  The case $a=0$ is trivial, and if $a\le -1$, then $-a\ge 1$ and what we have proved so far shows that
$(-a)^p\equiv (-a)\pmod p$.  But $(-a)^p=-a^p$ and hence $a^p\equiv a\pmod
p$.  We remark here that we do not believe the above proof of Theorem~\ref{fermat} is better than the standard group theoretic proof.  
Indeed, the ideas in it are seldom encountered in a first semester undergraduate level course in algebra.  However, the non-zero elements of $\Z_{p^j}$ do not form a group under multiplication if $j>1$, so that the standard argument does not generalize to the case of a power of a prime.  On the other hand, our method immediately gives a  proof of this case as well.

\begin{theorem} If $a\ge 1$ is any integer and $p$ is a prime, then
$a^{p^j}\equiv a \pmod p$ for all $j\ge 1$. \end{theorem}

{\bf Proof.} We note that a cyclic group $G$ of order $p^j$ acts on the set $A^{p^j}$ cyclically, and there are still precisely $a$ fixed points.  This time, however, there are $a^{p^j}$ total elements in the set. \qed

\section{The Euler Phi Function}

If $n\ge 1$ is an integer, we denote by $\phi(n)$ the number of elements in the set $\{1,\dots, n\}$ that are relatively prime to $n$.  The function $\phi$ is called the {\it Euler phi function}.  In this section, we use group actions to give a proof of the following well known result. 

\begin{theorem}\label{second}
If $n\ge 1$ is an integer, then $\displaystyle \sum_{d|n}\phi(d)=n$.
\end{theorem}
In contrast to Fermat's little theorem, this result is not typically discussed in algebra classes.  The usual proof given in number theory exploits the fact that the function $\phi$ is {\it multiplicative} (that is, satisfies $\phi(ab)=\phi(a)\phi(b)$ when ever $a$ and $b$ are relatively prime integers).  Our method here will employ Burnside's theorem as well as the structure of the lattice of subgroups of a finite cyclic group.  If a finite group $G$ acts on a finite set $S$, then for each $g\in G$, we let $S^g=\{s\in S : gs=s\}$ denote the set of elements in $S$ left fixed by $g$. If $r$ denotes the number of orbits in $S$ under the action of $G$, Burnside's theorem states that 
\begin{equation}
\label{burnside}
r\cdot |G|=\sum_{g\in G} |S^g|.
\end{equation}
We refer the reader to \cite{Fraleigh} for an excellent account of the details.  To establish Theorem~\ref{second}, we consider the problem of counting the number of distinguishable ways of coloring the edges of a regular $n$-gon ($n\ge 3$) with $q$ colors ($q\ge 1$).  The dihedral group $D_n$ of order $2n$ has a natural action on a regular $n$-gon as the group of symmetries.  Under this action, two colorings of the $n$-gon are indistinguishable if and only if they belong to the same orbit under the action.  Therefore the solution to our counting problem is the number of distinct orbits under this action.  

For notation, we let $D_n=\langle a,b | a^n=1, b^2=1, ba=a^{-1}b\rangle$.  We refer to an element of the cyclic subgroup $\langle a\rangle$ as a {\it rotation} and an element of the coset $b\langle a\rangle$ as a {\it flip}.  To use Burnside's theorem, we must compute $|S^g|$ for all $g\in D_n$ where $S$ is the set of all $q^n$ possible colorings of the $n$-gon.  If $g$ is a flip and $n$ is odd, then the line of reflection for $g$ must pass through a vertex of the $n$-gon and the midpoint of the edge opposite this vertex.  If a coloring $s$ is fixed under $g$, this opposite edge may be colored any one of the $q$ colors, but the remaining $(n-1)/2$ edges must be colored the same as their image under $g$.  Therefore there are $qq^{(n-1)/2}=q^{(n+1)/2}$ colorings fixed by $g$.  It follows that if $n$ is odd, then 
\[\sum_{g\in b\langle a\rangle} |S^g|=n\left (q^{(n+1)/2}\right).\]
If $n$ is even, then there are $q^{n/2}$ colorings fixed by $g$ if $g$ is a flip in a line through opposite vertices and there are $q^{(n+2)/2}$ colorings fixed by $g$ if $g$ is a flip in a line through the midpoints of opposite edges. Since there are exactly $n/2$ of each of these types of flips, if $n$ is even we have
\[\sum_{g\in b\langle a\rangle} |S^g|=\frac{n}{2}\left (q^{n/2}\right)+\frac{n}{2}q^{(n+2)/2}=\frac{n}{2}q^{n/2}(q+1).\]
Now we turn our attention toward the rotations.  For every positive divisor $d$ of $n=|\langle a\rangle|$, there is a unique subgroup of $\langle a\rangle$ of order $d$, and this subgroup has precisely $\phi(d)$ generators.  For each of these generators $g$, if we choose an edge of the $n$-gon, each of the images of this edge under the $d$ distinct powers of $g$ must be colored the same color if $g$ leaves the coloring fixed.  Therefore there are $q^{n/d}$ colorings left fixed by each of the $\phi(d)$ elements of order $d$ and hence we have
\[\sum_{g\in \langle a\rangle} |S^g|=\sum_{d|n}\phi(d)q^{n/d}.\]
Combining this with our results for the flips and Burnside's theorem, we have shown that the number of orbits in $S$ under the action of $D_n$ is given by 
\begin{equation*}
r=\left\{\begin{array}{ll} \frac{1}{2n}\left
(nq^{(n+1)/2}+{\displaystyle\sum_{d|n}\phi(d)q^{n/d}}\right ) & \mbox {if
$n$ is odd,}\\ \frac{1}{2n}\left
(\frac{n}{2}q^{n/2}(q+1)+{\displaystyle\sum_{d|n}\phi(d)q^{n/d}}\right ) &
\mbox {if $n$ is even.}\\ \end{array}\right. \end{equation*}
Now, Theorem~\ref{second} is easily verified for $n=1$ and $n=2$.  If $n\ge 3$, then setting $q=1$ above and noting that we must have $r=1$,  we see (using $n$ even or odd) we have
\[n+\sum_{d|n}\phi(d)=2n,\]
which completes the proof.

\bibliography{references}

\end{document}